\documentclass[10pt,11pt]{article}
\usepackage{amssymb}
\usepackage{amsmath}
\usepackage{amsfonts}

\evensidemargin=\oddsidemargin
\newdimen\dummy

\dummy=\oddsidemargin 
\input{tcilatex}

\begin{document}

\title{\textsf{Quasi-maximum likelihood estimation of periodic }$\mathsf{%
GARCH}$\textsf{\ processes}}
\author{\textsc{Abdelhakim Aknouche}$^{\text{*}}$ \and \textsc{Abdelouahab
Bibi}$^{\text{**}}$ \\
*Facult\'{e} de Math\'{e}matiques, Universit\'{e} U.S.T.H.B.\\
BP 32 El Alia, 16111, Bab ezzouar, Algiers, Algeria\\
E-mail: aknouche\_ab@yahoo.com\\
**D\'{e}partement de Math\'{e}matiques, Universit\'{e} Mentouri de
Constantine, Algeria\\
E-mail: abd.bibi@gmail.com\\
}
\date{}
\maketitle

\begin{abstract}
This paper establishes the strong consistency and asymptotic normality of
the quasi-maximum likelihood estimator ($QMLE$) for a $GARCH$ process with
periodically time-varying parameters. We first give a necessary and
sufficient condition for the existence of a strictly periodically stationary
solution for the periodic $GARCH$ ($P$-$GARCH$) equation. As a result, it is
shown that the moment of some positive order of the $P$-$GARCH$ solution is
finite, under which we prove the strong consistency and asymptotic normality 
$\left( CAN\right) $ of the $QMLE$ without any condition on the moments of
the underlying process.



\textit{Proposed running head}: $QML$ estimation of the periodic $GARCH$
\end{abstract}

\section{Introduction}

Periodic $GARCH$ ($P$-$GARCH$) processes introduced by Bollerslev and
Ghysels $(1996)$ have proved useful and appropriate for modeling many time
series encountered in practice, which are characterized by a stochastic
conditional variance with periodic dynamics (see also Franses and Paap $2000$%
; $2004$; Ghysels and Osborn, $2001$ and the references therein). As for the
standard $GARCH$ model, the $P$-$GARCH$ process may be seen as a non
Gaussian white noise but whose conditional variance follows a linear
periodic $ARMA$ ($PARMA$) dynamics in terms of the squared process, implying
a periodic structure for the underlying process itself. In contrast with
periodic $ARMA$ models which may be written as an equivalent vector $ARMA$
form (Tiao and Grupe, $1980$), there is no direct correspondence between the 
$P$-$GARCH$ model and the multivariate $GARCH$ one. In fact, any $P$-$GARCH$
process can be written as only a weak multivariate $GARCH$ model (see Drost
and Nijman $\left( 1993\right) $ for the definition of weak $GARCH$) meaning
that the study of the $P$-$GARCH$ may not be trivially deduced from the
existing multivariate $GARCH$ theory and thus it constitutes a useful and
interesting task.

Since their introduction, $P$-$GARCH$ models have been fairly considered in
the literature. Bollerslev and Ghysels $(1996)$ have studied the
second-order periodic stationarity of the $P$-$GARCH$ model and derived the
corresponding $QMLE$ but without studying its asymptotic properties. They
have successfully applied the $P$-$GARCH$ modeling for certain real data.
Franses and Paap $(2000)$ have considered the $P$-$GARCH$ process as an
innovation of a more general periodic $ARMA$ model with $P$-$GARCH$ error.
Important applications of such processes may be found in Ghysels and Osborn $%
(2001)$ and Franses and Paap $(2004)$. On the other hand, some probabilistic
properties such as strict\ periodic stationarity, existence of higher order
moment and geometric ergodicity have been studied recently by Bibi and
Aknouche $(2006)$. Apart from the mentioned works, it seems that there is no
result concerning asymptotic inference about such models. However this
problem has been exhaustively studied recently for the standard $GARCH$
case. Indeed, a considerable amount of research has been executed for
studying the asymptotic properties of the $QMLE$ for $GARCH$ processes. This
research, including contributions by authors such as Lee and Hansen $(1994)$%
, Lumsdaine $(1996)$, Boussama $(2000)$, Berkes et al. $(2003)$, Ling and
McAleer $(2003)$, Berkes and Horvath $(2004)$, and Francq and Zako\"{\i}an $%
(2004)$, is aimed at establishing consistency and asymptotic normality of
the $QMLE$\ for $GARCH$ processes with weak conditions on the moment as this
makes an undesirable restriction on the parameter space. Francq and Zako%
\"{\i}an $(2004)$'s result seems to be the weaker one obtained at present.

In the spirit of the latter work, the main goal of this paper is to
establish the strong consistency and asymptotic normality of the $QMLE$ for
the $P$-$GARCH$ process. We first show that the sufficient condition for
strict periodic stationarity of the $P$-$GARCH$ process given by Bibi and
Aknouche $(2006)$ is also necessary. As a consequence, it will be shown that
the moment of some positive order of the periodic stationary solution\ is
finite. This result will be exploited in order to obtain strong consistency
and asymptotic normality for the $QMLE$ irrespective of any moment
requirements for the underlying process.

The rest of this paper is organized as follows. Section 2 gives a necessary
and sufficient condition for the strict periodic stationarity for the $P$-$%
GARCH$ process and some important consequences thereof. The main results
concerning strong consistency and asymptotic normality of the $QMLE$ are
found in section 3 while the proofs of the main results are left in the
appendix.

Through the rest of this paper we need the use of the following notations: 
{\normalsize the symbols "$\leadsto "$ and "$a.s.$\ stand for convergence in
law and almost sure convergence respectively. A} matrix $A$ of order $%
m\times n$ is denoted by $A_{m\times n}$ and the $k\times k$ identity matrix
by $I_{k\times k}$. $\rho \left( A\right) $ denotes the spectral radius,
i.e. the maximum modulus of the eigenvalues of squared matrix $A$,$\
A^{\prime }$ stands for the transpose of $A$ and the vectorial relation $A>B$
means that each element of $A$ is greater than the corresponding element of $%
B$.

\section{Strict periodic stationarity of the $P$-$GARCH$ and some of its
important consequences}

This paper studies the periodic $GARCH$ process $\left( y_{t},t\in \mathbb{Z}%
\right) $ of orders $p$ and $q$ and period $S\geq 1$, that is a solution to
the stochastic difference equation%
\begin{equation}
\left\{ 
\begin{array}{l}
y_{t}=\sqrt{h_{t}}\eta _{t} \\ 
h_{t}=\omega _{t}+\sum\limits_{i=1}^{q}\alpha
_{t,i}y_{t-i}^{2}+\sum\limits_{j=1}^{p}\beta _{t,j}h_{t-j}%
\end{array}%
,\text{ }t\in \mathbb{Z}=\{0,\pm 1,\pm 2,\dots \},\right.  \tag{$2.1$}
\end{equation}%
where $\left( \eta _{t},t\in \mathbb{Z}\right) $ is a sequence of
independent and identically distributed $(i.i.d.)$ random variables defined
on a probability space $\left( \Omega ,\mathcal{A},P\right) $ such that $%
E\left( \eta _{t}\right) =0,$ $E\left( \eta _{t}^{2}\right) =1$ and $\eta
_{k}$ is independent of $y_{t}$ for all $k>t$. The parameters $\omega _{t}$, 
$\alpha _{t,i},$ and $\beta _{t,j}$ are periodic in $t$ with period $S$
(i.e., $\omega _{t+kS}=\omega _{t}$, $\alpha _{t+kS,i}=\alpha _{t,i}$ and $%
\beta _{t+kS,j}=\beta _{t,j}$) such that $\omega _{t}>0,$ $\alpha _{t,i}\geq
0$ and $\beta _{t,j}\geq 0,$ $i=1,...,q,$ $j=1,...,p$, for all $k,t\in 
\mathbb{Z}$, so that by setting $t=v+Sn$, model $(2.1)$ may be equivalently
written as%
\begin{equation}
\left\{ 
\begin{array}{l}
y_{v+Sn}=\sqrt{h_{v+Sn}}\eta _{v+Sn} \\ 
h_{v+Sn}=\omega _{v}+\sum\limits_{i=1}^{q}\alpha
_{v,i}y_{v+Sn-i}^{2}+\sum\limits_{j=1}^{p}\beta _{v,j}h_{v+Sn-j}%
\end{array}%
,t\in \mathbb{Z},\right.  \tag{$2.2$}
\end{equation}%
highlighting thus the periodicity in the model. In the difference equation $%
\left( 2.2\right) $ $y_{v+Sn}$ may refer to $\left( y_{t}\right) $ during
the $v-th$ `season' of the year $n$ and $\omega _{v},\alpha _{v,1},...,$ $%
\alpha _{v,q}$ together with $\beta _{v,1},...,\beta _{v,p}$ are the model
coefficients at season $v$. The orders $p$ and $q$ can also be considered
periodic in time, say $p_{t}$ and $q_{t}$ in order to reduce the number of
parameters to be estimated. However, in this case a constrained estimation
procedure, which is not the objective of this paper, would be needed.

A standard approach for studying the structure of an $ARCH$-type model is to
write it as a random coefficient representation. Indeed, it is easy to write
model $(2.1)$ in term of the squared process as follows%
\begin{equation}
\left\{ 
\begin{array}{l}
y_{t}^{2}=\omega _{t}\eta _{t}^{2}+\sum\limits_{i=1}^{q}\alpha _{t,i}\eta
_{t}^{2}y_{t-i}^{2}+\sum\limits_{j=1}^{p}\beta _{t,j}\eta _{t}^{2}h_{t-j} \\ 
h_{t}=\omega _{t}+\sum\limits_{i=1}^{q}\alpha
_{t,i}y_{t-i}^{2}+\sum\limits_{j=1}^{p}\beta _{t,j}h_{t-j},%
\end{array}%
\right. \text{ }t\in \mathbb{Z},  \tag{$2.3$}
\end{equation}%
which is ready to be cast in a stochastic recurrence equation with random
coefficients.

Indeed, defining the $\left( p+q\right) $-random vectors $%
Y_{t}=(y_{t}^{2},...,y_{t-q+1}^{2},h_{t},...,h_{t-p+1})^{\prime }$ and

$B_{t}=\left( \omega _{t}\eta _{t}^{2},0_{\left( q-1\right) \times
1}^{\prime },\omega _{t},0_{\left( p-1\right) \times 1}^{\prime }\right)
^{\prime }$, ($0_{m\times n}$ stands for the null matrix of order $m\times n$%
) together with the $(p+q)\times (p+q)$ random matrix $A_{t}$ given by%
\begin{equation*}
A_{t}=\left( 
\begin{array}{ll}
\left. \alpha _{t,1}\eta _{t}^{2}\right. \left. \alpha _{t,2}\eta
_{t}^{2}\right. \ldots \left. \alpha _{t,q-1}\eta _{t}^{2}\right. \left.
\alpha _{t,q}\eta _{t}^{2}\right. & \left. \beta _{t,1}\eta _{t}^{2}\right.
\left. \beta _{t,2}\eta _{t}^{2}\right. \ldots \left. \beta _{t,p-1}\eta
_{t}^{2}\right. \left. \beta _{t,p}\eta _{t}^{2}\right. \\ 
\left. \text{ \ \ \ \ \ }\right. I_{\left( q-1\right) \times \left(
q-1\right) }\text{ }\left. \text{ \ \ \ \ \ \ \ \ }\right. 0_{\left(
q-1\right) \times 1} & \left. \text{ \ \ \ \ \ }\right. 0_{\left( q-1\right)
\times p} \\ 
\left. \alpha _{t,1}\right. \left. \text{ \ \ }\alpha _{t,2}\text{ \ \ }%
\right. \left. \ldots \right. \left. \alpha _{t,q-1}\right. \text{ \ \ \ \ }%
\left. \alpha _{t,q}\right. & \left. \beta _{t,1}\text{ \ }\right. \left. 
\text{ \ }\beta _{t,2}\text{ \ \ }\right. \left. \ldots \right. \left. \beta
_{t,p-1}\right. \left. \text{ \ \ }\beta _{t,p}\right. \\ 
\left. \text{ \ \ \ \ \ \ \ \ \ }\right. 0_{\left( p-1\right) \times q} & 
\left. \text{ \ \ \ \ \ \ \ \ \ \ }\right. I_{\left( p-1\right) \times
\left( p-1\right) }\left. \text{\ \ \ \ \ \ }\right. 0_{\left( p-1\right)
\times 1}%
\end{array}%
\right) ,
\end{equation*}%
one can rewrite model $\left( 2.3\right) $\ in the following generalized $AR$
model%
\begin{equation}
Y_{t}=A_{t}Y_{t-1}+B_{t},\text{ }t\in \mathbb{Z},  \tag{$2.4$}
\end{equation}%
which differs from the standard formulation studied by Bougerol and Picard ($%
1992a$, $1992b$) in that the coefficients $(A_{t},B_{t})$ are rather \textit{%
independent and periodically distributed} ($i.p.d.$). Moreover, their
expectations $E\left\{ A_{t}\right\} $ and $E\left\{ B_{t}\right\} $ defined
element-wise are $S$-periodic in time.

Now, we focus on the stochastic recurrence equation $(2.4)$ and the
existence of a \textit{strictly periodically stationary} (henceforth $s.p.s.$%
) solution thereof. Recall that the process $\left( Y_{t},\text{ }t\in 
\mathbb{Z}\right) $ is said to be $s.p.s.$ (with period $S\geq 1$) if the
distribution of $\left( Y_{t_{1}},Y_{t_{2}},...,Y_{t_{n}}\right) ^{\prime }$
is the same as that of $\left(
Y_{t_{1}+Sh},Y_{t_{2}+Sh},...,Y_{t_{n}+Sh}\right) ^{\prime }$ for all $n\geq
1$ and $h$, $t_{1},t_{2},...,t_{n}\in \mathbb{Z}$. Furthermore, it is called 
\textit{periodically ergodic} (cycloergodic, cf, Boyles and Gardner, $1983$)
if\ for all Borel set $B$ and all integer $m$%
\begin{equation*}
\frac{1}{n}\sum_{t=1}^{n}\mathbb{I}_{B}\left(
Y_{v+St},Y_{v+1+St},...,Y_{v+m+St}\right) \longrightarrow P\left( \left(
Y_{v},Y_{v+1},...,Y_{v+m}\right) \in B\right) ,\text{ }a.s.\text{ as }%
n\rightarrow \infty ,
\end{equation*}%
for all $1\leq v\leq S$, where $\mathbb{I}_{B}(.)$ denotes the indicator
function of the set $B$. It follows immediately that a standard ergodic
process is periodically ergodic with period $S=1$. As for the stationary
case (Billingsley, 1996 Theorem 36.4), periodic ergodicity is closed under
certain transformations. In particular if $\left( \varepsilon _{t},\text{ }%
t\in \mathbb{Z}\right) $ is $s.p.s.$ and periodically ergodic and if $\left(
Y_{t},\text{ }t\in \mathbb{Z}\right) $ is given by $Y_{t}=f\left(
...,\varepsilon _{t},\varepsilon _{t+1},...\right) $ where $f$ is a
measurable function from $\mathbb{R}^{\infty }$ to $\mathbb{R}$, then $%
\left( Y_{t},\text{ }t\in \mathbb{Z}\right) $ is also $s.p.s.$ and
periodically ergodic (Aknouche and Guerbyenne, $2007$). A periodic analog of
the ergodic theorem for stationary sequences can be stated as follows. If $%
\left( Y_{t},\text{ }t\in \mathbb{Z}\right) $ is $s.p.s.$ and periodically
ergodic and if $f$ is a measurable function from $\mathbb{R}^{\infty }$ to $%
\mathbb{R}$ such that $E\left\{ f\left( ...,Y_{t-1},Y_{t},Y_{t+1},...\right)
\right\} <\infty $ then for all $1\leq v\leq S$,%
\begin{equation*}
\frac{1}{n}\sum_{t=1}^{n}f\left( ...,Y_{v+S\left( t-1\right)
},Y_{v+St},Y_{v+S\left( t+1\right) },...\right) \rightarrow E\left\{ f\left(
...,Y_{v+S\left( t-1\right) },Y_{v+St},Y_{v+S\left( t+1\right) },...\right)
\right\} ,\text{ }a.s.\text{ as }n\rightarrow \infty .
\end{equation*}%
It is well known that the periodic process theory can be translated
immediately into the stationarity process one by an appropriate
transformation. Since the seminal paper by Gladyshev $(1961)$, we know that
a periodically stationary process $\left( Y_{t},\text{ }t\in \mathbb{Z}%
\right) $ is equivalent to a vector-valued stationary process $\left( 
\underline{Y}_{n},\text{ }n\in \mathbb{Z}\right) $ where 
\begin{equation*}
\underline{Y}_{n}=\left( Y_{1+nS},Y_{2+nS},...,Y_{S+nS}\right) ^{\prime },
\end{equation*}%
meaning that\textbf{\ }any property about periodic processes translates at
once into a corresponding result about stationary processes. In our case,
the process $\left( \underline{Y}_{n},n\in \mathbb{Z}\right) $ satisfies
nonuniquely the stochastic difference equation%
\begin{equation}
\underline{Y}_{n}=\underline{A}_{n}\underline{Y}_{n-1}+\underline{B}_{n},%
\text{ \ \ \ \ }n\in \mathbb{Z},  \tag{$2.5$}
\end{equation}%
where $\underline{A}_{n}$ and $\underline{B}_{n}$ are defined by blocks as%
\begin{equation*}
\underline{A}_{n}=\left( 
\begin{array}{cccc}
0_{r\times r} & \ldots & 0_{r\times r} & A_{1+nS} \\ 
0_{r\times r} & \ldots & 0_{r\times r} & A_{2+nS}A_{1+nS} \\ 
\vdots & \vdots & \vdots & \vdots \\ 
0_{r\times r} & \ldots & 0_{r\times r} & \prod\limits_{v=0}^{S-1}A_{nS+S-j}%
\end{array}%
\right) _{rS\times rS},\underline{B}_{n}=\left( 
\begin{array}{c}
B_{1+nS} \\ 
A_{2+nS}B_{1+nS}+B_{2+nS} \\ 
\vdots \\ 
\sum\limits_{k=1}^{S}\left\{ \prod\limits_{v=0}^{S-k-1}A_{S-v+nS}\right\}
B_{k+nS}%
\end{array}%
\right) _{rS\times 1},
\end{equation*}%
with $r=p+q$ and where, as usual, empty products are set equal to one.%
\textbf{\ }Therefore, the solution process $(2.4)$ is $s.p.s.$ if and only
if the corresponding solution of $(2.5)$ is strictly stationary. Similarly,
the process given by $(2.4)$ is periodically ergodic if and only if the
corresponding process solution of $(2.5)$ is ergodic. So, we can substitute
the study of the properties of model $(2.4)$ to those of model $(2.5)$ if
they are much easy to obtain. In the sequel we only use model $(2.4)$ as it
seems to be simpler.

As for strict stationarity, the primordial tool in studying strict periodic
stationarity is the top Lyapunov exponent for $i.p.d.$ random matrices,
properties of which can be found in Aknouche and Guerbyenne $(2007)$. Let $%
\left\Vert .\right\Vert $\ be an arbitrary operator norm in $\mathcal{M}^{r}$%
, the space of real square matrices of dimension $r$. Then, the top Lyapunov
exponent associated with the\ $i.p.d.$\ sequence of matrices $A:=\left(
A_{t},\text{ }t\in \mathbb{Z}\right) $ is defined when $\sum_{v=1}^{S}E%
\left( \log ^{+}\left\Vert A_{v}\right\Vert \right) <\infty $,\ by%
\begin{equation}
\gamma ^{S}\left( A\right) =\inf_{n\in \mathbb{N}^{\ast }}\frac{1}{n}%
E\left\{ \log \left\Vert A_{nS}A_{nS-1}...A_{1}\right\Vert \right\} , 
\tag{$2.6$}
\end{equation}%
where $\log ^{+}(x)=\max (\log \left( x\right) ,0)$. For $S=1$, $(2.6)$
reduces to the definition of the top Lyapunov exponent for $i.i.d.$ matrices
(e.g. Bougerol et Picard, $1992a$). It is clear that $\gamma ^{S}\left(
.\right) $ inherits the properties of the standard top Lyapunov exponent. In
particular, the following inequality $\gamma ^{S}\left( A\right) \leq
\sum_{v=1}^{S}E\left( \log \left\Vert A_{v}\right\Vert \right) $ holds with
equality for the scalar case $r=1$. Furthermore, when the $\left( A_{t},t\in 
\mathbb{Z}\right) $ is a sequence of $S$-periodic nonrandom matrices, then
from $(2.6)$ 
\begin{equation*}
\gamma ^{S}\left( A\right) =\log \lim\limits_{n\rightarrow \infty }\left(
\left\Vert \left( \prod_{v=0}^{S-1}A_{S-v}\right) ^{n}\right\Vert \right)
^{1/n}\overset{def}{=}\log \rho \left( \prod_{v=1}^{S}A_{v}\right) .
\end{equation*}%
On the other hand, by the $i.p.d.$ property of $\left( A_{t},\text{ }t\in 
\mathbb{Z}\right) $, the sequence$\left( \prod_{v=0}^{S-1}A_{nS-v},n\in 
\mathbb{Z}\right) $\ is $i.i.d.$ to which applying the subadditive ergodic
theorem (cf., Kingman $1973$) or the result of Furstemberg and Kesten $%
(1960) $ it follows that $a.s.$%
\begin{equation}
\gamma ^{S}\left( A\right) =\lim_{n\rightarrow \infty }\frac{1}{n}\log
\left\Vert A_{nS}A_{nS-1}...A_{1}\right\Vert .  \tag{$2.7$}
\end{equation}%
As noted below, properties of model $(2.4)$ can be obtained from those of
model $(2.5)$ which is a generalized $AR$ with nonnegative $i.i.d.$
coefficients for which a fairly exhaustive theory exists (see Bougerol and
Picard $1992a$, $1992b$ and the references therein). Since $\left( \eta
_{t},t\in \mathbb{Z}\right) $ is $i.i.d.$ with finite variance and the
parameter are $S$-periodic, the sequence $\left( \left( A_{t},B_{t}\right)
,t\in \mathbb{Z}\right) $ is $s.p.s.$ and periodically ergodic and then the
random pair $\left( \left( \underline{A}_{t},\underline{B}_{t}\right) ,t\in 
\mathbb{Z}\right) $ is stationary and ergodic. Moreover, because $%
\sum_{v=1}^{S}E\left( \log ^{+}\left\Vert A_{v}\right\Vert \right) \leq
\sum_{v=1}^{S}E\left( \left\Vert A_{v}\right\Vert \right) $ and $%
\sum_{v=1}^{S}E\left( \log ^{+}\left\Vert B_{v}\right\Vert \right) \leq
\sum_{v=1}^{S}E\left( \left\Vert B_{v}\right\Vert \right) $, then $%
\sum_{v=1}^{S}E\left( \log ^{+}\left\Vert A_{v}\right\Vert \right) $ and $%
\sum_{v=1}^{S}E\left( \log ^{+}\left\Vert B_{v}\right\Vert \right) $
together with $E\left( \log ^{+}\left\Vert \underline{A}_{0}\right\Vert
\right) $ and $E\left( \log ^{+}\left\Vert \underline{B}_{0}\right\Vert
\right) $ are finite. Thus from Theorem 2.3 of Bougerol and Picard $(1992b)$
a necessary and sufficient condition for model $(2.5)$ to possess a \textit{%
nonanticipative} (future-independent) strictly stationary solution is that
the top Lyapunov exponent associated with $\underline{A}=\left( \underline{A}%
_{t},t\in \mathbb{Z}\right) $, that is $\gamma \left( \underline{A}\right)
:=\inf_{t\in \mathbb{N}^{\ast }}\frac{1}{t}E\left\{ \left( \log \left\Vert 
\underline{A}_{t}\underline{A}_{t-1}...\underline{A}_{1}\right\Vert \right)
\right\} $ is strictly negative.

It must be noted that the latter result is not suitable since it gives a
condition about the transformed non-periodic model $(2.5)$, not the original
periodic one $(2.4)$, which is the objective of this paper. Fortunately,
from the relation between $\left( \underline{A}_{t},t\in \mathbb{Z}\right) $
and $\left( A_{t},t\in \mathbb{Z}\right) $, it is easy to see, taking a
multiplicative norm, that $\gamma \left( \underline{A}\right) \leq \gamma
^{S}\left( A\right) $, so that a sufficient condition for model $(2.4)$ to
possess a nonanticipative $s.p.s.$ solution is that $\gamma ^{S}\left(
A\right) <0$, which is the result of Theorem 2 of Bibi and Aknouche $(2006)$%
. Nevertheless, using model $(2.5)$ it seems to be difficult to show that
the latter condition is also necessary. Using directly model $(2.4)$, the
following theorem shows that this condition is also necessary.

\noindent \textbf{Theorem 2.1} \textit{Model }$(2.4)$\textit{\ admits a
nonanticipative }$s.p.s.$\textit{\ solution given by}%
\begin{equation}
Y_{t}=\sum\limits_{k=1}^{\infty
}\prod\limits_{j=0}^{k-1}A_{t-j}B_{t-k}+B_{t},\text{ }t\in \mathbb{Z}, 
\tag{$2.8$}
\end{equation}%
\textit{\ if and only if the top Lyapunov exponent }$\gamma ^{S}\left(
A\right) $\textit{\ given by }$(2.6)$\textit{\ is strictly negative, where
the above series converges }$a.s$\textit{. for all }$t\in \mathbb{Z}$\textit{%
. Moreover, the solution is unique and periodically ergodic.}

We now give a necessary condition for the strict stationarity, which
coincide for the case $S=1$ (non periodic) to corollary 2.3 of Bougerol and
Picard $(1992b)$ and which we will be needed in studying asymptotics for the 
$QMLE$.

\noindent \textbf{Corollary 2.2} \textit{If the }$P$\textit{-}$GARCH$\textit{%
\ model }$(2.1)$\textit{\ possesses an }$s.p.s.$\textit{\ solution then}%
\begin{equation*}
\rho \left( \prod\limits_{v=1}^{S}\mathbf{\beta }_{v}\right) <1,
\end{equation*}%
\textit{where }$\mathbf{\beta }_{t}$\textit{\ is the submatrix of }$A_{t}$%
\textit{\ defined by }%
\begin{equation*}
\mathbf{\beta }_{t}=\left( 
\begin{array}{cc}
\beta _{t,1}\beta _{t,2}\ldots \beta _{t,p-1} & \beta _{t,p} \\ 
I_{\left( p-1\right) \times \left( p-1\right) } & 0_{\left( p-1\right)
\times 1}%
\end{array}%
\right) .
\end{equation*}

We now turn to an important consequence of strict periodic stationarity and
the corresponding log-moment condition of Theorem 2.1, that is the existence
of a moment of some positive order. A similar result for the non periodic $%
GARCH$ case was proved by Nelson $\left( 1990\right) $ for the $GARCH(1,1)$
and Berkes et al. $\left( 2003\right) $ (Lemma $2.3$) for the general $GARCH$
case.

\noindent \textbf{Theorem 2.3} \textit{If }$\gamma ^{S}\left( A\right) <0$%
\textit{\ then there is }$\delta >0$\textit{\ such that}%
\begin{equation}
E\left( h_{t}^{\delta }\right) <\infty \text{ \ \ \ and \ \ }E\left(
y_{t}^{2\delta }\right) <\infty \text{.}  \tag{$2.9$}
\end{equation}

\section{Asymptotic properties of the quasi-maximum likelihood estimator}

Let $\theta =(\theta _{1}^{\prime },\theta _{2}^{\prime },...,\theta
_{S}^{\prime })^{\prime }$, where $\theta _{i}=(\omega _{i},\alpha
_{i,1},...,\alpha _{i,q},\beta _{i,1},...,\beta _{i,p})^{\prime }$, be the
parameter vector which is supposed to belong to a compact parameter space $%
\Theta \subset \left( \left] 0,+\infty \right[ \times \lbrack 0,+\infty
\lbrack ^{(p+q)}\right) ^{S}$. The true parameter value is unknown and is
denoted by $\theta ^{0}=(\theta _{1}^{0\prime },\theta _{2}^{0\prime
},...,\theta _{S}^{0\prime })^{\prime }$ with $\theta _{i}^{0}=(\omega
_{i}^{0},\alpha _{i,1}^{0},...,\alpha _{i,q}^{0},\beta _{i,1}^{0},...,\beta
_{i,p}^{0})^{\prime }$. Consider a series $\left\{
y_{1,}y_{2},...,y_{T},T=NS\right\} $ generated from model $(2.1)$ with true
parameter $\theta 
{{}^\circ}%
\in \Theta $, i.e. from the model%
\begin{equation}
\left\{ 
\begin{array}{l}
y_{v+Sn}=\sqrt{h_{v+Sn}}\eta _{v+Sn} \\ 
h_{v+Sn}=\omega _{v}^{0}+\sum\limits_{i=1}^{q}\alpha
_{v,i}^{0}y_{v+Sn-i}^{2}+\sum\limits_{j=1}^{p}\beta _{v,j}^{0}h_{v+Sn-j}%
\end{array}%
\text{ \ }t\in \mathbb{Z},\right.  \tag{$3.1$}
\end{equation}%
where the sample size $T$ is supposed without loss of generality multiple of
the period $S$. The Gaussian log-likelihood function of $\theta \in \Theta $
conditional on initial values $y_{0},...,y_{1-q},\widetilde{h}_{0},...,%
\widetilde{h}_{1-p}$ is given (ignoring a constant) by%
\begin{equation}
\widetilde{L}_{NS}\left( \theta \right) =-\frac{1}{NS}\sum_{k=1}^{N}%
\sum_{s=1}^{S}\widetilde{l}_{s+kS}\left( \theta \right) ,  \tag{$3.2$}
\end{equation}%
with%
\begin{equation}
\widetilde{l}_{t}\left( \theta \right) =\frac{y_{t}^{2}}{\widetilde{h}_{t}}%
+\log \widetilde{h}_{t},  \tag{$3.3$}
\end{equation}%
where $\widetilde{h}_{t}$ is solution to the conditional model%
\begin{equation}
\left\{ 
\begin{array}{l}
y_{v+Sn}=\sqrt{\widetilde{h}_{v+Sn}}\eta _{v+Sn} \\ 
\widetilde{h}_{v+Sn}=\omega _{v}+\sum\limits_{i=1}^{q}\alpha
_{v,i}y_{v+Sn-i}^{2}+\sum\limits_{j=1}^{p}\beta _{v,j}\widetilde{h}_{v+Sn-j}%
\end{array}%
,\text{ \ }t\geq 1,\right.  \tag{$3.4$}
\end{equation}%
with given initial values $y_{0},...,y_{1-q},\widetilde{h}_{0},...,%
\widetilde{h}_{1-p}.$ These values may be chosen, taking into account the $S$%
-periodicity of the distribution of $\left( y_{t}\right) $, as an
approximation of the non conditional variance. For instance for the first
order $P$-$GARCH$ model the non conditional variance is given by%
\begin{equation*}
E\left( y_{t}^{2}\right) =E\left( h_{t}\right) =\frac{\omega
_{t}+\sum\limits_{j=1}^{S-1}\prod\limits_{i=0}^{j-1}\left( \alpha
_{t-i}+\beta _{t-i}\right) \omega _{t-j}}{1-\prod\limits_{i=0}^{S-1}\left(
\alpha _{t-i}+\beta _{t-i}\right) },
\end{equation*}%
which may be negative for some values of the strict periodic stationarity
domain. Thus, taking $\alpha _{t}=\beta _{t}=0$, we have the initial values%
\begin{equation}
\left. 
\begin{array}{c}
y_{0}^{2}=\widetilde{h}_{0}=\omega _{0}, \\ 
y_{-1}^{2}=\widetilde{h}_{-1}=\omega _{-1}^{2}, \\ 
\vdots \\ 
y_{1-d}=\widetilde{h}_{1-d}=\omega _{-d}^{2},%
\end{array}%
\right.  \tag{$3.5a$}
\end{equation}%
which can be considered the same for the general $P$-$GARCH$ model. Another
choice of $y_{0},...,y_{1-q},\widetilde{h}_{0},...,$ $\widetilde{h}_{1-p}$ is%
\begin{equation}
\left. 
\begin{array}{c}
y_{0}^{2}=\widetilde{h}_{0}=y_{S}^{2}, \\ 
y_{-1}^{2}=\widetilde{h}_{-1}=y_{\left[ S-1\right] }^{2}, \\ 
\vdots \\ 
y_{1-d}=\widetilde{h}_{1-d}=y_{\left[ S-d\right] }^{2},%
\end{array}%
\right.  \tag{$3.5b$}
\end{equation}%
where $d=\max (p,q)$ and $\left[ k\right] =k$ if $k\geq 1$ and $\left[ k%
\right] =lS+k$ otherwise, $l$ being the lower integer such that $lS+k\geq 1$%
. Obviously, these choices have only a practical value and do not affect the
asymptotic properties.

The $QMLE$ of $\theta 
{{}^\circ}%
$, denoted by $\widehat{\theta }_{NS}$, is the maximizer of $\widetilde{L}%
_{NS}\left( \theta \right) $ on $\Theta $ and then the minimizer of%
\begin{equation*}
\frac{1}{NS}\sum_{k=1}^{N}\sum_{s=1}^{S}\widetilde{l}_{s+kS}\left( \theta
\right) .
\end{equation*}%
Let $\gamma ^{S}\left( A^{0}\right) $ denote the top Lyapunov exponent
associated with $\left( A_{t}^{0},t\in \mathbb{Z}\right) $ where $A_{t}^{0}$
is just $A_{t}$ defined in $(2.4)$ with $\theta ^{0}$ in place of $\theta $.
To study the strong consistency of $\widehat{\theta }_{NS}$ consider the
following assumptions.

\noindent \textbf{A1}: $\gamma ^{S}\left( A^{0}\right) <0$ and $\forall
\theta \in \Theta ,$ $\rho \left( \prod\limits_{v=1}^{S}\mathbf{\beta }%
_{v}\right) <1$.

\noindent \textbf{A2}: The polynomials $\alpha _{v}^{0}\left( z\right)
=\sum\limits_{j=1}^{q}\alpha _{v,j}^{0}z^{j}$ and $\beta _{v}^{0}\left(
z\right) =1-\sum\limits_{j=1}^{p}\beta _{v,j}^{0}z^{j}$\ have no common
roots,\ $\alpha _{v}^{0}\left( 1\right) \neq 0$ and $\alpha _{v,q}^{0}+\beta
_{v,p}^{0}\neq 0$ for all $1\leq v\leq S$.

\noindent \textbf{A3}: $\left( \eta _{t}\right) $ is non degenerate.

As seen in Theorem 2.3, assumption \textbf{A1} ensures the existence of a
finite moment, for the solution of $(3.1)$, which is the key for proving
consistency and asymptotic normality irrespective of any moment condition.
The condition $\rho \left( \prod\limits_{v=1}^{S}\mathbf{\beta }_{v}\right)
<1$ is imposed for any $\theta \in \Theta $ in order to obtain $h_{t}\left(
\theta \right) $ as a causal solution of $\left\{ y_{t},y_{t-1},...\right\} $
while \textbf{A2} is made to guarantee the identifiability of the parameter.
Similar conditions for the $GARCH$ case have been considered by Francq and
Zako\"{\i}an $(2004)$.

Given a realization $\left\{ y_{1},...,y_{NS}\right\} $,\ $\widetilde{l}%
_{t}\left( \theta \right) $ can be approximated for $1\leq t\leq NS$ by 
\begin{equation}
l_{t}\left( \theta \right) =\dfrac{\varepsilon _{t}^{2}}{h_{t}\left( \theta
\right) }+\log h_{t}\left( \theta \right) ,  \tag{$3.6$}
\end{equation}%
where $h_{t}\left( \theta \right) $ is the $s$.$p$.$s$. solution of%
\begin{equation}
h_{v+Sn}=\omega _{v}+\sum\limits_{i=1}^{q}\alpha
_{v,i}y_{v+Sn-i}^{2}+\sum\limits_{j=1}^{p}\beta _{v,j}h_{v+Sn-j},\text{ }%
t\in \mathbb{Z},  \tag{$3.7$}
\end{equation}%
for $\theta \in \Theta $ with $h_{t}=h_{t}\left( \theta 
{{}^\circ}%
\right) $. Let 
\begin{equation*}
L_{NS}\left( \theta \right) =-\frac{1}{NS}\sum_{k=1}^{N}%
\sum_{s=1}^{S}l_{s+kS}\left( \theta \right) .
\end{equation*}%
The following theorem establishes the strong consistency of $\widehat{\theta 
}_{NS}$.

\noindent \textbf{Theorem 3.1}

\noindent \textit{Under }\textbf{A1}-\textbf{A3}\textit{\ }$\widehat{\theta }%
_{NS}$\textit{\ is strongly consistent in the sense that }$\widehat{\theta }%
_{NS}\rightarrow \theta 
{{}^\circ}%
$ $a.s.$ as $N\rightarrow \infty .$

\noindent \textbf{Remarks}

\noindent 1) From assumptions \textbf{A1-A3}, it is clear that the above
result remains true for the particular periodic $ARCH$ ($P$-$ARCH$) case,
i.e. when $\beta _{t,j}=0$.

\noindent 2) In order that \textbf{A1} holds for, for example, the
particular $P$-$GARCH(1,1)$, we may chose the parameter space $\Theta $ as a
compact of the form $\Theta =\left( [\epsilon ,1/\epsilon ]\times \lbrack
0,1/\epsilon ]\times \lbrack 0,1-\epsilon ]\right) ^{S}$ where $\epsilon >0$
is so small that the parameter $\theta ^{0}=(\omega _{1}^{0},\alpha
_{1,1}^{0},\beta _{1,1}^{0},\omega _{2}^{0},\alpha _{2,1}^{0},\beta
_{2,1}^{0}...,\omega _{S}^{0},\alpha _{S,1}^{0},\beta _{S,1}^{0})^{\prime }$
belongs to $\Theta $. The fact that $\beta _{v,1}^{0}$ must be lower than 1
for all $v$ is so restrictive and can be weakened so that $%
\prod\limits_{v=1}^{S}\beta _{v,1}^{0}<1$ which is just assumption \textbf{A1%
} for the particular $P$-$GARCH(1,1)$.

Now we turn to the asymptotic normality of $\widehat{\theta }_{NS}$.
Consider the following additional assumptions.

\noindent \textbf{A4}: $\theta 
{{}^\circ}%
$ is in the interior of $\Theta $.

\noindent \textbf{A5}: $E\left( \eta _{t}^{4}\right) <\infty $.

Assumption \textbf{A4} is standard and allows to validate the first order
condition on the maximizer of the log-likelihood while \textbf{A5} is
necessary for the existence of the limiting covariance matrix of the $QMLE$.

The following result establishes the asymptotic normality of $\widehat{%
\theta }_{NS}$.

\noindent \textbf{Theorem 3.2}

\noindent \textit{Under \textbf{A1-A5}\ we have}%
\begin{equation}
\sqrt{NS}\left( \widehat{\theta }_{NS}-\theta 
{{}^\circ}%
\right) \leadsto N\left( 0,\left( E\left( \eta _{t}^{4}\right) -1\right)
J^{-1}\right) \text{, as }N\rightarrow \infty ,  \tag{$3.9$}
\end{equation}%
\textit{where the matrix }$J$\textit{\ given by}%
\begin{equation}
J:=\left[ \sum_{v=1}^{S}E_{\theta 
{{}^\circ}%
}\left( \frac{\partial ^{2}l_{v}\left( \theta 
{{}^\circ}%
\right) }{\partial \theta \partial \theta ^{\prime }}\right) \right]
=\sum_{v=1}^{S}E_{\theta 
{{}^\circ}%
}\left( \frac{1}{h_{v}^{2}\left( \theta 
{{}^\circ}%
\right) }\frac{\partial h_{v}\left( \theta 
{{}^\circ}%
\right) }{\partial \theta }\frac{\partial h_{v}\left( \theta 
{{}^\circ}%
\right) }{\partial \theta ^{\prime }}\right) ,  \tag{$3.10$}
\end{equation}%
\textit{is block diagonal}.

Let us now apply the forgoing results to the first order $P$-$ARCH$ process
given by%
\begin{equation*}
y_{v+Sn}=\eta _{v+Sn}\sqrt{\omega _{v}^{0}+\alpha _{v}^{0}y_{v+Sn-1}^{2}},
\end{equation*}%
where $\omega _{v}^{0}>0$ and $\alpha _{v}^{0}>0$. It is easily seen that
the strict periodic stationarity condition for the $P$-$ARCH(1)$ reduces to%
\begin{equation*}
0\leq \prod\limits_{v=1}^{S}\alpha _{v}^{0}<\exp \left( -E\left( \log \left(
\eta _{v}^{2}\right) \right) \right) :=a,
\end{equation*}%
under which, supposing that $\theta ^{0}=(\omega _{1}^{0},\alpha
_{1,1}^{0},\omega _{2}^{0},\alpha _{2,1}^{0},...,\omega _{S}^{0},\alpha
_{S,1}^{0})^{\prime }$ belongs to a compact $\Theta $ of the form%
\begin{equation}
\Theta =\left( [\epsilon ,1/\epsilon ]\times \lbrack 0,a^{1/S}-1]\right)
^{S},  \tag{$3.11$}
\end{equation}%
\ the $QMLE$ is then by Theorem 3.1 strongly consistent.

Moreover, if we assume that $\theta ^{0}$ is in the interior of $\Theta $
given by $(3.11)$, then from Theorem 3.2 the $QMLE$ is also asymptotically
Gaussian and its limiting distribution is given by $(3.9)$ and $(3.10)$
which in the $P$-$GARCH(1,1)$ case reduces to%
\begin{equation*}
J=\left( 
\begin{array}{cccc}
J_{1} & 0_{2\times 2} & \cdots & 0_{2\times 2} \\ 
0_{2\times 2} & J_{2} & \cdots & \vdots \\ 
\vdots & \vdots & \ddots & 0_{2\times 2} \\ 
0_{2\times 2} & \cdots & 0_{2\times 2} & J_{S}%
\end{array}%
\right) ,
\end{equation*}%
where 
\begin{equation*}
J_{v}=E\left( 
\begin{array}{cc}
\frac{1}{\left( \omega _{v}+\alpha _{v,1}\varepsilon _{v-1}^{2}\right) ^{2}}
& \frac{\varepsilon _{v-1}^{2}}{\left( \omega _{v}+\alpha _{v,1}\varepsilon
_{v-1}^{2}\right) ^{2}} \\ 
\frac{\varepsilon _{v-1}^{2}}{\left( \omega _{v}+\alpha _{v,1}\varepsilon
_{v-1}^{2}\right) ^{2}} & \frac{\varepsilon _{v-1}^{4}}{\left( \omega
_{v}+\alpha _{v,1}\varepsilon _{v-1}^{2}\right) ^{2}}%
\end{array}%
\right) ,\text{ }1\leq v\leq S.
\end{equation*}

\section{Appendix. Proofs}

\subsection{Appendix A. Proofs of Section 2's results}

\subsubsection{Proof of Theorem 2.1}

We only give a proof of the necessity part of the theorem as the sufficiency
one was established earlier by Bibi and Aknouche $(2006)$. Suppose that
model $(2.4)$ admits a nonanticipative $s.p.s.$ solution $\left( Y_{t},t\in 
\mathbb{Z}\right) $. Then iterating $(2.4)$ $m$ -times for some $m>1$ and
some $v\in \left\{ 1,...,S\right\} $ gives%
\begin{equation*}
Y_{v}=\dsum\limits_{j=0}^{m}\prod\limits_{i=0}^{j-1}A_{v-i}B_{v-j}+\prod%
\limits_{j=0}^{m+1}A_{v-j}Y_{v-m-1},
\end{equation*}%
and exploiting the non-negativity of the coefficients of $A_{t}$ and those
of $Y_{v}$ it follows that for all $m>1$%
\begin{equation*}
Y_{v}\geq \dsum\limits_{j=0}^{m}\prod\limits_{i=0}^{j-1}A_{v-i}B_{v-j},\text{
\ }a.s.
\end{equation*}%
implying that the series $\dsum\limits_{j=0}^{\infty
}\prod\limits_{i=0}^{j-1}A_{v-i}B_{v-j}$ converge $a.s.$ and therefore%
\begin{equation}
\prod\limits_{i=0}^{j-1}A_{v-i}B_{v-j}\rightarrow 0,\text{ \ }a.s.\text{ as }%
j\rightarrow \infty ,  \tag{$A.1$}
\end{equation}%
from which we want to show that%
\begin{equation}
\prod\limits_{i=0}^{j-1}A_{v-i}\rightarrow 0,\text{ \ }a.s.\text{ as }%
j\rightarrow \infty .  \tag{$A.2$}
\end{equation}%
This holds whenever 
\begin{equation}
\lim\limits_{j\rightarrow \infty }\prod\limits_{i=0}^{j-1}A_{v-i}e_{k}=0,a.s.%
\text{ for all }k\in \left\{ 1,...,r\right\} ,  \tag{$A.3$}
\end{equation}%
where $\left( e_{k}\right) $ is the canonical basis of $\mathbb{R}^{r}$.
Observing that $B_{v-j}=\omega _{v-j}\eta _{v-j}^{2}e_{1}+\omega
_{v-j}e_{q+1}$, $(A.1)$ implies that%
\begin{equation}
\omega _{v-j}\eta _{v-j}^{2}\prod\limits_{i=0}^{j-1}A_{v-i}e_{1}\rightarrow 0%
\text{, }a.s.\text{ as }j\rightarrow \infty \text{, and }\omega
_{v-j}\prod\limits_{i=0}^{j-1}A_{v-i}e_{q+1}\rightarrow 0\text{, }a.s.\text{
as }j\rightarrow \infty ,  \tag{$A.4$}
\end{equation}%
and since $\omega _{v-j}>0$, this shows $(A.3)$ for $k=q+1$. Moreover, since
for all $k=1,...,p$%
\begin{equation*}
A_{v-j+1}e_{q+k}=\beta _{v-j+1,k}\eta _{v-j+1}^{2}e_{1}+\beta
_{v-j+1,k}e_{q+1}+e_{q+k+1},
\end{equation*}%
then for $k=1$ we have almost surely 
\begin{equation*}
0=\lim\limits_{j\rightarrow \infty
}\prod\limits_{i=0}^{j-1}A_{v-i}e_{q+1}\geq \lim\limits_{j\rightarrow \infty
}\prod\limits_{i=0}^{j-1}A_{v-i}e_{q+2}\geq 0,
\end{equation*}%
showing $(A.3)$ for $k=q+2$ and then recursively for $k=q+j$, $j=1...p$. On
the other hand, since 
\begin{equation*}
A_{v-j+1}e_{q}=\alpha _{v-j+1,q}\eta _{v-j+1}^{2}e_{1}+\alpha
_{v-j+1,q}e_{q+1},
\end{equation*}%
then from $(A.4),$ follows $(A3)$ for $k=q$. Finally, from the identity 
\begin{equation*}
A_{v-j+1}e_{k}=\alpha _{v-j+1,k}\eta _{v-j+1}^{2}e_{1}+\alpha
_{v-j+1,k}e_{q+1}+e_{k+1},\text{ \ \ }k=1,...,q-1,
\end{equation*}%
follows $(A.3)$ for the rest values of $k$ using a backward recursion. This
show $(A.2)$, which implies that any subsequence of the sequence\ $%
U_{j}:=\left( \prod\limits_{i=0}^{j-1}A_{v-i}\right) _{j}$ converge $a.s.$
to zero as $j\rightarrow \infty .$ In particular, 
\begin{equation}
U_{Sj}=\prod\limits_{i=0}^{Sj-1}A_{v-i}=\prod\limits_{i=0}^{j-1}\left(
\prod\limits_{k=0}^{S-1}A_{v-iS-k}\right) \rightarrow 0,\text{ \ }a.s.\text{
as }j\rightarrow \infty .  \tag{$A.5$}
\end{equation}%
Now, since the sequence $\left( \prod\limits_{k=0}^{S-1}A_{v-iS-k},i\in 
\mathbb{Z}\right) $ is $i.i.d.$ for all $v$, then from $(A.5)$ and Lemma 2.1
of Bougerol and Picard $(1992b)$ it follows that the top Lyapunov exponent
associated with the $i.i.d.$ sequence $\left(
\prod\limits_{k=0}^{S-1}A_{v-iS-k},i\in \mathbb{Z}\right) $, which is
exactly $\gamma ^{S}\left( A\right) $, is strictly negative. This completes
the proof.

\subsubsection{Proof of Corollary 2.2}

\noindent As the $\left( A_{t}\right) $ is nonnegative, the corresponding
top Lyapunov exponent $\gamma ^{S}\left( A\right) $ is greater that the top
Lyapunov exponent corresponding to nonrandom sequence of matrices, say $%
\mathbf{\beta =}\left( \mathbf{\beta }_{t},t\in \mathbb{Z}\right) $ obtained
when one set the elements of the first $q$ lines and those of the first $q$
columns to be equal to zero. In other words, 
\begin{equation*}
\gamma ^{S}\left( A\right) \geq \gamma ^{S}\left( \mathbf{\beta }\right)
:=\log \rho \left( \prod_{v=1}^{S}\mathbf{\beta }_{v}\right) .
\end{equation*}%
Thus, if the model has an $s.p.s.$ solution which is equivalent to the
condition $\gamma ^{S}\left( A\right) <0$, the conclusion of the Corollary
is true.

\subsubsection{Proof of Theorem 2.3}

\noindent The proof is similar to that of Lemma 2.3 of Berkes et al. $(2003)$%
. First we have to show that if $\gamma ^{S}\left( A\right) <0$ then there
is $\delta >0$ and $n_{0}$ such that 
\begin{equation}
E\left( \left\Vert A_{n_{0}S}A_{n_{0}S-1}...A_{1}\right\Vert ^{\delta
}\right) <1.  \tag{$A.6$}
\end{equation}%
Since $\gamma ^{S}\left( A\right) =\inf_{n\in \mathbb{N}^{\ast }}\left\{ 
\frac{1}{n}E\left( \log \left\Vert A_{nS}A_{nS-1}...A_{1}\right\Vert \right)
\right\} $ is strictly negative, there is a positive integer $n_{0}$ such
that%
\begin{equation*}
E\left( \log \left\Vert A_{n_{0}S}A_{n_{0}S-1}...A_{1}\right\Vert \right) <0.
\end{equation*}%
On the other hand,\ working with a multiplicative norm and by the property $%
i.p.d.$ of the sequence $\left( A_{t}\right) $ we have%
\begin{eqnarray*}
E\left( \left\Vert A_{n_{0}S}A_{n_{0}S-1}...A_{1}\right\Vert \right)
&=&\left\Vert E\left( A_{n_{0}S}A_{n_{0}S-1}...A_{1}\right) \right\Vert \\
&=&\left\Vert E\left( A_{S}A_{S-1}...A_{1}\right) ^{n_{0}}\right\Vert \\
&\leq &\left\Vert E\left( A_{S}A_{S-1}...A_{1}\right) \right\Vert
^{n_{0}}<\infty .
\end{eqnarray*}%
Let $f\left( t\right) =E\left( \left\Vert
A_{n_{0}S}A_{n_{0}S-1}...A_{1}\right\Vert ^{t}\right) $. Since $f^{\prime
}\left( 0\right) =E\left( \log \left\Vert
A_{n_{0}S}A_{n_{0}S-1}...A_{1}\right\Vert \right) <0,$ $f\left( t\right) $
decrease in a neighborhood of $0$ and since $f\left( 0\right) =1,$ it
follows that there exists $0<\delta <1$ such that $\left( A.6\right) $ holds.

Now from Theorem 2.1 we have for some $v\in \left\{ 1,...,S\right\} $%
\begin{equation*}
\left\Vert Y_{v}\right\Vert \leq \sum\limits_{k=1}^{\infty }\left\Vert
\prod\limits_{j=0}^{k-1}A_{v-j}\right\Vert \left\Vert B_{v-k}\right\Vert
+\left\Vert B_{v}\right\Vert ,
\end{equation*}%
and since $0<\delta <1$ then%
\begin{equation*}
\left\Vert Y_{v}\right\Vert ^{\delta }\leq \sum\limits_{k=1}^{\infty
}\left\Vert \prod\limits_{j=0}^{k-1}A_{v-j}\right\Vert ^{\delta }\left\Vert
B_{v-k}\right\Vert ^{\delta }+\left\Vert B_{v}\right\Vert ^{\delta },
\end{equation*}%
which by the independence of $A_{v-j}$ and $B_{v-k}$ for $j<k$ implies 
\begin{eqnarray*}
E\left\Vert Y_{v}\right\Vert ^{\delta } &\leq &\sum\limits_{k=1}^{\infty
}E\left( \left\Vert \prod\limits_{j=0}^{k-1}A_{v-j}\right\Vert ^{\delta
}\right) E\left( \left\Vert B_{v-k}\right\Vert ^{\delta }\right) +E\left(
\left\Vert B_{v}\right\Vert ^{\delta }\right) \\
&\leq &B\left( \delta \right) \sum\limits_{k=1}^{\infty }E\left( \left\Vert
\prod\limits_{j=0}^{k-1}A_{v-j}\right\Vert ^{\delta }\right) +E\left(
\left\Vert B_{v}\right\Vert ^{\delta }\right) ,
\end{eqnarray*}%
where $B\left( \delta \right) =\max\limits_{0\leq v\leq S-1}E\left(
\left\Vert B_{v-k}\right\Vert ^{\delta }\right) $. Using $(A.6)$ there exist 
$a_{v}>0$ and $0<b_{v}<1$ such that%
\begin{equation*}
E\left( \left\Vert \prod\limits_{j=0}^{k-1}A_{v-j}\right\Vert ^{\delta
}\right) \leq a_{v}b_{v}^{k}\leq ab^{k},
\end{equation*}%
where $ab^{k}=\max\limits_{0\leq v\leq S-1}\left\{ a_{v}b_{v}^{k}\right\} $.
This prove that $E\left\Vert Y_{v}\right\Vert ^{\delta }<\infty $ and then $%
(2.9)$.

\subsection{Appendix B. Proofs of Theorems 3.1-3.2}

The proofs of Theorems 3.1 and 3.2 are by now standard and follow from
similar arguments used in showing the strong consistency and asymptotic
normality of the $QMLE$ for standard $GARCH$ models (Francq and Zako\"{\i}%
an, $2004$). The main aim here is to reveal the basic assumptions and to
quantify the asymptotic distribution of the $QMLE$ for the $P$-$GARCH$.
Since there are several similarities between the standard $GARCH$ and the $P$%
-$GARCH$, certain steps of the proof for the $P$-$GARCH$ case are similar in
spirit to that of the standard $GARCH$ one. Thus, we give details of proof
only when it seems pertinent and refer to Francq and Zako\"{\i}an $(2004)$
for the details.

\subsubsection{Proof of Theorem 3.1}

Theorem 3.1 will be proved by showing several lemmas. Lemma B.1 establishes
the uniform asymptotic forgetting of initial values, Lemma B.2 ensures
identifiability of the parameter, Lemma B.3 shows the finiteness of the
limiting criterion $\sum_{s=1}^{S}E_{\theta 
{{}^\circ}%
}\left( l_{s}\left( \theta \right) \right) $ and that this one is uniquely
minimized at $\theta 
{{}^\circ}%
$, while Lemma B.4 uses the compactness of $\Theta $\ and a periodically
ergodic argument to conclude the strong convergence.

\noindent \textbf{Lemma B.1} \textit{Under \textbf{A1} we have}%
\begin{equation*}
\underset{\theta \in \Theta }{\sup }\left\vert L_{NS}\left( \theta \right) -%
\widetilde{L}_{NS}\left( \theta \right) \right\vert \rightarrow 0,a.s.\text{
as }N\rightarrow \infty ,
\end{equation*}

\noindent \textbf{Proof} Rewriting $(3.7)$ in a vectorial form as follows%
\begin{equation}
\underline{h}_{t}=\mathbf{\beta }_{t}\underline{h}_{t-1}+\underline{\mathbf{%
\alpha }}_{t},\text{ }t\in \mathbb{Z},  \tag{$A.7$}
\end{equation}%
where $\underline{h}_{t}=\left( h_{t},h_{t-1},...,h_{t-p+1}\right) ^{\prime
} $ and $\underline{\mathbf{\alpha }}_{t}=\left( \omega
_{v}+\sum\limits_{i=1}^{q}\alpha _{v,i}y_{v+Sn-i}^{2},0,...,0\right)
_{1\times p}^{\prime }$. From Corollary 2.2, the assumption $\rho \left(
\prod\limits_{v=1}^{S}\mathbf{\beta }_{v}\right) <1$ of \textbf{A1} and the
compactness of $\Theta $ we have%
\begin{equation}
\sup_{\theta \in \Theta }\rho \left( \prod\limits_{v=1}^{S}\mathbf{\beta }%
_{v}\right) <1,  \tag{$A.8$}
\end{equation}%
which implies, by iterating $(A.7)$, that%
\begin{equation*}
\underline{h}_{t}=\sum_{k=0}^{t-1}\prod\limits_{i=0}^{k-1}\mathbf{\beta }%
_{t-i}\underline{\mathbf{\alpha }}_{t-k}+\prod\limits_{i=0}^{t}\mathbf{\beta 
}_{t-i}\underline{h}_{0},\text{ }t\in \mathbb{Z}.
\end{equation*}%
If we denote by $\widetilde{\underline{h}}_{t}$ and $\underline{\widetilde{%
\mathbf{\alpha }}}_{t}$ the vectors obtained from $\underline{h}_{t}$ and $%
\underline{\mathbf{\alpha }}_{t}$,\ respectively, while replacing $h_{t-j}$\
by $\widetilde{h}_{t-j}$\ with initial values given by $(3.5)$, then we have%
\begin{equation*}
\widetilde{\underline{h}}_{t}=\sum_{k=0}^{t-q-1}\prod\limits_{i=0}^{k-1}%
\mathbf{\beta }_{t-i}\underline{\mathbf{\alpha }}_{t-k}+\sum_{k=t-q}^{t-1}%
\prod\limits_{i=0}^{k-1}\mathbf{\beta }_{t-i}\underline{\widetilde{\mathbf{%
\alpha }}}_{t-k}+\prod\limits_{i=0}^{t}\mathbf{\beta }_{t-i}\widetilde{%
\underline{h}}_{0},
\end{equation*}%
and from $(A.8)$ it follows that%
\begin{eqnarray*}
\sup_{\theta \in \Theta }\left\Vert \underline{h}_{t}-\widetilde{\underline{h%
}}_{t}\right\Vert &=&\sup_{\theta \in \Theta }\left\Vert
\sum_{k=t-q}^{t-1}\prod\limits_{i=0}^{k-1}\mathbf{\beta }_{t-i}\left( 
\underline{\mathbf{\alpha }}_{t-k}-\underline{\widetilde{\mathbf{\alpha }}}%
_{t-k}\right) +\prod\limits_{i=0}^{t}\mathbf{\beta }_{t-i}\left( \underline{h%
}_{0}-\widetilde{\underline{h}}_{0}\right) \right\Vert \\
&\leq &M\rho ^{t},
\end{eqnarray*}%
from which, we get%
\begin{eqnarray*}
\underset{\theta \in \Theta }{\sup }\left\vert L_{NS}\left( \theta \right) -%
\widetilde{L}_{NS}\left( \theta \right) \right\vert &\leq &\frac{1}{N}%
\sum_{t=1}^{NS}\sup_{\theta \in \Theta }\left[ \left\vert \frac{\widetilde{h}%
_{t}-h_{t}}{\widetilde{h}_{t}h_{t}}\right\vert y_{t}^{2}+\left\vert \log
\left( \frac{h_{t}}{\widetilde{h}_{t}}\right) \right\vert \right] \\
&\leq &\left( \max_{1\leq v\leq S}\sup_{\theta \in \Theta }\left( \omega
_{v}^{-2}\right) \right) \frac{M}{N}\sum_{t=1}^{NS}\rho ^{t}y_{t}^{2}+\left(
\max_{1\leq v\leq S}\sup_{\theta \in \Theta }\left( \omega _{v}^{-1}\right)
\right) \frac{M}{N}\sum_{t=1}^{NS}\rho ^{t},
\end{eqnarray*}%
where the inequality $\left\vert \log \frac{y}{x}\right\vert \leq \frac{%
\left\vert y-x\right\vert }{\min \left( y,x\right) }$ for positive $x$ and $%
y $ has been used. The existence of a moment for $y_{t}$ (cf, Corollary 2.2)
implies, by the Borel-Cantelli lemma, that $\rho ^{t}y_{t}^{2}\rightarrow 0$ 
$a.s.$ so that the conclusion of the lemma follows from the Toeplitz lemma.

\noindent \textbf{Lemma B.2} \textit{Under \textbf{A1-A3} there is }$t\in 
\mathbb{Z}$\textit{\ such that if }$h_{t}\left( \theta \right) =h_{t}\left(
\theta 
{{}^\circ}%
\right) $\textit{\ }$a$\textit{.}$s$\textit{. then }$\theta =\theta 
{{}^\circ}%
$.

\noindent \textbf{Proof} From the assumption $\rho \left(
\prod\limits_{v=1}^{S}\mathbf{\beta }_{v}\right) <1$ of \textbf{A1} the
polynomials $\left( \beta _{v}^{0}\left( L\right) \right) _{1\leq v\leq S}$ (%
$L$ being the backward shift operator) are invertible. Suppose that $%
h_{t}\left( \theta \right) =h_{t}\left( \theta 
{{}^\circ}%
\right) $\textit{\ }$a$\textit{.}$s$\textit{.} for some $t$, then using $%
(3.7)$ we have%
\begin{equation*}
\left( \frac{\alpha _{v}\left( L\right) }{\beta _{v}\left( L\right) }-\frac{%
\alpha _{v}^{0}\left( L\right) }{\beta _{v}^{0}\left( L\right) }\right)
y_{v+St}^{2}=\left( \frac{\omega _{v}^{0}}{\beta _{v}^{0}\left( 1\right) }-%
\frac{\omega _{v}}{\beta _{v}\left( 1\right) }\right) ,\text{ for all\ }%
1\leq v\leq S.
\end{equation*}%
If $\dfrac{\alpha _{v}\left( L\right) }{\beta _{v}\left( L\right) }\neq 
\dfrac{\alpha _{v}^{0}\left( L\right) }{\beta _{v}^{0}\left( L\right) }$\
for some\ $1\leq v\leq S$\ then there exists a deterministic periodic
time-varying combination of $y_{v+St-j}^{2}$, $j\geq 1$. This contradict the
fact that $y_{v+St}^{2}=E\left( y_{v+St}^{2}/y_{v+St-1}^{2},...\right)
+h_{v+St}\left( \eta _{v+St}^{2}-1\right) $, since by \textbf{A3} $\left(
\eta _{t},t\in \mathbb{Z}\right) $ is non degenerate. Therefore,%
\begin{equation*}
\frac{\alpha _{v}\left( z\right) }{\beta _{v}\left( z\right) }=\frac{\alpha
_{v}^{0}\left( z\right) }{\beta _{v}^{0}\left( z\right) }\text{ }\forall
\left\vert z\right\vert \leq 1\text{\ and }\frac{\omega _{v}^{0}}{\beta
_{v}^{0}\left( 1\right) }-\frac{\omega _{v}}{\beta _{v}\left( 1\right) },%
\text{\ for all\ }1\leq v\leq S,
\end{equation*}%
which by the assumption \textbf{A2} of absence of common roots implies that $%
\alpha _{v}\left( z\right) =\alpha _{v}^{0}\left( z\right) $, $\beta
_{v}\left( z\right) =\beta _{v}^{0}\left( z\right) $ and $\omega _{v}=\omega
_{v}^{0}$ for all\ $1\leq v\leq S$, proving the lemma.

\noindent \textbf{Lemma B.3} \textit{Under}\textbf{\ }\textit{\textbf{A1}}%
\begin{equation*}
\sum_{v=1}^{S}E_{\theta 
{{}^\circ}%
}\left( l_{v}\left( \theta ^{0}\right) \right) <\infty ,
\end{equation*}%
\textit{\ and }$\sum_{v=1}^{S}E_{\theta 
{{}^\circ}%
}\left( l_{v}\left( \theta \right) \right) $\textit{\ is minimum at }$\theta
=\theta 
{{}^\circ}%
$.

\noindent \textbf{Proof} By Theorem 2.3 we have%
\begin{eqnarray*}
\sum_{v=1}^{S}E_{\theta 
{{}^\circ}%
}\left( \log h_{v}\left( \theta ^{0}\right) \right)
&=&\sum_{v=1}^{S}E_{\theta 
{{}^\circ}%
}\frac{1}{\delta }\left( \log h_{v}\left( \theta ^{0}\right) ^{\delta
}\right) \\
&\leq &\frac{1}{\delta }\sum_{v=1}^{S}\log E_{\theta 
{{}^\circ}%
}\left( h_{v}\left( \theta ^{0}\right) ^{\delta }\right) <\infty ,
\end{eqnarray*}%
from which it follows that%
\begin{eqnarray*}
\sum_{v=1}^{S}E_{\theta 
{{}^\circ}%
}\left( l_{v}\left( \theta ^{0}\right) \right) &=&\sum_{v=1}^{S}E_{\theta 
{{}^\circ}%
}\left[ \frac{h_{v}\left( \theta ^{0}\right) \eta _{v}}{h_{v}\left( \theta
^{0}\right) }+\log h_{v}\left( \theta ^{0}\right) \right] \\
&=&S+\sum_{v=1}^{S}E_{\theta 
{{}^\circ}%
}\left( \log h_{v}\left( \theta ^{0}\right) \right) <\infty .
\end{eqnarray*}%
Finally,%
\begin{eqnarray}
\sum_{v=1}^{S}E_{\theta 
{{}^\circ}%
}\left( l_{v}\left( \theta \right) \right) -\sum_{v=1}^{S}E_{\theta 
{{}^\circ}%
}\left( l_{v}\left( \theta ^{0}\right) \right) &=&\sum_{v=1}^{S}E_{\theta 
{{}^\circ}%
}\left[ \log \frac{h_{v}\left( \theta \right) }{h_{v}\left( \theta
^{0}\right) }+\frac{h_{v}\left( \theta ^{0}\right) }{h_{v}\left( \theta
\right) }-1\right]  \notag \\
&\geq &\sum_{v=1}^{S}E_{\theta 
{{}^\circ}%
}\left[ \log \frac{h_{v}\left( \theta \right) }{h_{v}\left( \theta
^{0}\right) }+\log \frac{h_{v}\left( \theta ^{0}\right) }{h_{v}\left( \theta
\right) }\right] =0,  \TCItag{$A.9$}
\end{eqnarray}%
showing that the limit criterion is minimized at $\theta ^{0}$.

\noindent \textbf{Lemma B.4} \textit{Under}\textbf{\ }\textit{\textbf{A1}},
for all $\theta \neq \theta 
{{}^\circ}%
$ there is a neighborhood $\mathcal{V}\left( \theta \right) $ such that%
\begin{equation*}
\lim_{N\rightarrow \infty }\inf \inf_{\widetilde{\theta }\in \mathcal{V}%
\left( \theta \right) }\left( -\frac{1}{N}\widetilde{L}_{NS}\left( 
\widetilde{\theta }\right) \right) >\sum_{v=1}^{S}E_{\theta 
{{}^\circ}%
}\left( l_{v}\left( \theta ^{0}\right) \right) .
\end{equation*}%
\noindent \textbf{Proof} For all $\theta \in \Theta $ and all integer $k$,
let $\mathcal{V}_{k}\left( \theta \right) $ an open sphere of center $\theta 
$ and radius $1/k$. Using Lemma B.1 we have%
\begin{eqnarray*}
\lim_{N\rightarrow \infty }\inf \inf_{\widetilde{\theta }\in \mathcal{V}%
_{k}\left( \theta \right) \cap \Theta }\left( -\frac{1}{N}\widetilde{L}%
_{NS}\left( \widetilde{\theta }\right) \right) &\geq &%
\begin{array}{l}
\lim_{N\rightarrow \infty }\inf \inf_{\widetilde{\theta }\in \mathcal{V}%
_{k}\left( \theta \right) \cap \Theta }\left( -\frac{1}{N}L_{NS}\left( 
\widetilde{\theta }\right) \right) - \\ 
\text{ \ \ \ \ \ \ }\lim_{N\rightarrow \infty }\sup \frac{1}{N}\sup_{\theta
\in \Theta }\left\vert \widetilde{L}_{NS}\left( \theta \right) -L_{NS}\left(
\theta \right) \right\vert%
\end{array}
\\
&\geq &\lim_{N\rightarrow \infty }\inf \frac{1}{N}\sum_{n=0}^{N-1}%
\sum_{v=1}^{S}\inf_{\widetilde{\theta }\in \mathcal{V}_{k}\left( \theta
\right) \cap \Theta }l_{v+nS}\left( \widetilde{\theta }\right) .
\end{eqnarray*}%
Applying the ergodic theorem for the $i.i.d.$ sequence $\left(
\sum_{v=1}^{S}l_{v+nS}\left( \widetilde{\theta }\right) \right) _{n}$ with $%
E\left( \sum_{v=1}^{S}l_{v+nS}\left( \widetilde{\theta }\right) \right) \in 
\mathbb{R\cup }\left\{ \infty \right\} $ (cf, Billingsley $1995$, p. $284,$ $%
495$) it follows that%
\begin{equation*}
\lim_{N\rightarrow \infty }\inf \frac{1}{N}\sum_{n=0}^{N-1}\sum_{v=1}^{S}%
\inf_{\widetilde{\theta }\in \mathcal{V}_{k}\left( \theta \right) \cap
\Theta }l_{v+nS}\left( \widetilde{\theta }\right) =\sum_{v=1}^{S}E_{\theta 
{{}^\circ}%
}\left( \inf_{\widetilde{\theta }\in \mathcal{V}_{k}\left( \theta \right)
\cap \Theta }l_{v+nS}\left( \widetilde{\theta }\right) \right) ,
\end{equation*}%
and by the Beppo-Levi theorem (e.g. Billingsley, $1995,$ p. 219), we have%
\begin{equation*}
\sum_{v=1}^{S}E_{\theta 
{{}^\circ}%
}\left( \inf_{\widetilde{\theta }\in \mathcal{V}_{k}\left( \theta \right)
\cap \Theta }l_{v+nS}\left( \widetilde{\theta }\right) \right) \rightarrow
\sum_{v=1}^{S}E_{\theta 
{{}^\circ}%
}\left( l_{v}\left( \theta \right) \right) \text{ as }k\rightarrow \infty ,
\end{equation*}%
which by $(A.9)$ proves the lemma.

\noindent \textbf{Proof of Theorem 3.1} In view of Lemmas B.1-B.4, the proof
of the theorem is completed by using an argument of compactness of $\Theta $%
. First, for all neighborhood $\mathcal{V}\left( \theta ^{0}\right) $ of $%
\theta ^{0}$ we have%
\begin{eqnarray}
\lim_{N\rightarrow \infty }\sup \inf_{\widetilde{\theta }\in \mathcal{V}%
\left( \theta ^{0}\right) }\left( -\frac{1}{N}\widetilde{L}_{NS}\left( 
\widetilde{\theta }\right) \right) &\leq &\lim_{N\rightarrow \infty }\left( -%
\frac{1}{N}\widetilde{L}_{NS}\left( \theta ^{0}\right) \right)
=\lim_{N\rightarrow \infty }\left( -\frac{1}{N}L_{NS}\left( \theta
^{0}\right) \right)  \notag \\
&=&\sum_{v=1}^{S}E_{\theta 
{{}^\circ}%
}\left( l_{v}\left( \theta ^{0}\right) \right) .  \TCItag{$A.10$}
\end{eqnarray}%
The compact $\Theta $ is recovered by a union of a neighborhood $\mathcal{V}%
\left( \theta ^{0}\right) $ of $\theta ^{0}$ and the set of neighborhoods $%
\mathcal{V}\left( \theta \right) $, $\theta \in \Theta \backslash \mathcal{V}%
\left( \theta ^{0}\right) $, where $\mathcal{V}\left( \theta \right) $
fulfills Lemma B.4. Therefore, there exists a finite sub-covering of $\Theta 
$ by $\mathcal{V}\left( \theta _{0}\right) ,$ $\mathcal{V}\left( \theta
_{1}\right) ,...,\mathcal{V}\left( \theta _{k}\right) $ such that%
\begin{equation*}
\inf_{\theta \in \mathcal{\Theta }}\left( -\frac{1}{N}\widetilde{L}%
_{NS}\left( \widetilde{\theta }\right) \right) =\min_{i\in \left\{
1,2,...,k\right\} }\inf_{\widetilde{\theta }\in \Theta \cap \mathcal{V}%
\left( \theta _{i}\right) }\left( -\frac{1}{N}\widetilde{L}_{NS}\left( 
\widetilde{\theta }\right) \right) .
\end{equation*}%
From $(A.10)$ and Lemma B.4, the latter relation shows that $\widehat{\theta 
}_{NS}\in \mathcal{V}\left( \theta ^{0}\right) $ for $N$ sufficiently large,
which complete the proof of the theorem.

\subsubsection{Proof of Theorem 3.2}

The proof of Theorem 3.2 is based on a Taylor expansion of $\dfrac{\partial 
\widetilde{L}_{NS}\left( \theta \right) }{\partial \theta }$ at $\theta 
{{}^\circ}%
$, along the lines of Francq and Zako\"{\i}an $(2004)$, which is given by%
\begin{equation*}
0=\left( N\right) ^{-1/2}\sum_{t=1}^{NS}\dfrac{\partial \widetilde{l}_{t}}{%
\partial \theta }\left( \widehat{\theta }_{NS}\right) =\left( N\right)
^{-1/2}\sum_{t=1}^{NS}\dfrac{\partial \widetilde{l}_{t}}{\partial \theta }%
\left( \theta 
{{}^\circ}%
\right) +\left( \left( N\right) ^{-1}\sum_{t=1}^{NS}\dfrac{\partial ^{2}%
\widetilde{l}_{t}}{\partial \theta \partial \theta ^{\prime }}\left( 
\widetilde{\theta }\right) \right) \left( N\right) ^{1/2}\left( \widehat{%
\theta }_{NS}-\theta 
{{}^\circ}%
\right) ,
\end{equation*}%
where the coordinates of $\widetilde{\theta }$ are between the corresponding
entrees of $\widehat{\theta }_{NS}$ and those of $\theta 
{{}^\circ}%
$.

The proof is also based on several lemmas, Lemma B.5-B.9, which are aimed at
establishing: the integrability of the first derivatives of the limiting
criterion at $\theta 
{{}^\circ}%
$, the invertibility of $J$ and its relation with de first derivatives of
the limiting criterion, the uniform integrability of the third-order
derivatives of the limiting criterion, the asymptotic forgetting of starting
values for the derivatives, and a central limit theorem for martingales
differences together with a periodically ergodic theorem for the second
derivatives of the criterion.

\noindent \textbf{Lemma B.5} \textit{We have}%
\begin{equation*}
\sum_{s=1}^{S}E_{\theta 
{{}^\circ}%
}\left\vert \dfrac{\partial l_{s}\left( \theta 
{{}^\circ}%
\right) }{\partial \theta }\dfrac{\partial l_{s}\left( \theta 
{{}^\circ}%
\right) }{\partial \theta ^{\prime }}\right\vert <\infty \text{, and\ }%
\sum_{s=1}^{S}E_{\theta 
{{}^\circ}%
}\left\vert \dfrac{\partial ^{2}l_{s}\left( \theta 
{{}^\circ}%
\right) }{\partial \theta \partial \theta ^{\prime }}\right\vert <\infty .
\end{equation*}

\noindent \textbf{Lemma B.6} \textit{Under \textbf{A1-A5} }$J$\textit{\ is
invertible and}%
\begin{equation*}
\sum_{s=1}^{S}E_{\theta 
{{}^\circ}%
}\left\vert \dfrac{\partial l_{s}\left( \theta 
{{}^\circ}%
\right) }{\partial \theta }\right\vert ^{2}=\left( E\left( \eta
_{t}^{4}\right) -1\right) J.
\end{equation*}

\noindent \textbf{Lemma B.7} \textit{The following limit relations}%
\begin{equation*}
N^{-\frac{1}{2}}\sum_{k=1}^{N}\sum_{s=1}^{S}\left\vert \frac{\partial
l_{s+kS}\left( \theta 
{{}^\circ}%
\right) }{\partial \theta }-\frac{\partial \widetilde{l}_{s+kS}\left( \theta 
{{}^\circ}%
\right) }{\partial \theta }\right\vert =o_{p}\left( 1\right) ,
\end{equation*}%
\textit{and}%
\begin{equation*}
\underset{\theta \in v\left( \theta 
{{}^\circ}%
\right) }{\sup }N^{-\frac{1}{2}}\sum_{k=1}^{N}\sum_{s=1}^{S}\left\vert \frac{%
\partial l_{s+kS}\left( \theta 
{{}^\circ}%
\right) }{\partial \theta }-\frac{\partial \widetilde{l}_{s+kS}\left( \theta 
{{}^\circ}%
\right) }{\partial \theta }\right\vert =o_{p}\left( 1\right) ,
\end{equation*}%
\textit{\ are true}.

\noindent \textbf{Lemma B.8} \textit{There is a neighborhood }$V\left(
\theta ^{0}\right) $\textit{\ of }$\theta ^{0}$\textit{\ such that for all }$%
i,j,k\in \left\{ 1,...,S(p+q+1)\right\} $%
\begin{equation*}
\sum_{s=1}^{S}E_{\theta 
{{}^\circ}%
}\sup_{\theta \in \mathcal{V}\left( \theta ^{0}\right) }\left\vert \dfrac{%
\partial ^{3}l_{s}\left( \theta \right) }{\partial \theta _{i}\partial
\theta _{j}\partial \theta _{k}}\right\vert <\infty .
\end{equation*}

\noindent \textbf{Lemma B.9} \textit{The following limit results are true}%
\begin{equation*}
N^{-\frac{1}{2}}\sum_{k=1}^{N}\sum_{s=1}^{S}\frac{\partial l_{s+kS}}{%
\partial \theta }\left( \theta ^{0}\right) \leadsto N\left( 0,\left( E\left(
\eta _{t}^{4}\right) -1\right) J\right) \text{ and }N^{-1}\sum_{k=1}^{N}%
\sum_{s=1}^{S}\frac{\partial ^{2}l_{s+kS}}{\partial \theta \partial \theta
^{\prime }}\left( \widetilde{\theta }\right) \rightarrow J\text{, }a.s.
\end{equation*}%
The proofs are very similar to those of Francq and Zako\"{\i}an $(2004)$. It
suffices to replace the stationarity and ergodicity arguments by the
periodic stationarity and periodic ergodicity ones, respectively. For this
we omit the proofs which are quite lengthy.

\end{document}